\documentclass{amsart}

\newtheorem{theorem}{Theorem}[section]
\newtheorem{lemma}[theorem]{Lemma}
\newtheorem{corollary}[theorem]{Corollary}

\theoremstyle{definition}

\theoremstyle{remark}
\newtheorem{remark}[theorem]{Remark}

\numberwithin{equation}{section}

\def \N{I\!\!N}

\def \P{I\!\!P}
\def \T{T\!\!\!T}

\def \d{\, \nu_g}
\def \s{\int_{M}}
\def \dis {\displaystyle}
\def \la{\lambda_{1,p}}
\def \mf{|u|^{p-2}}
\def \dfp{\int_M |du|^p\, \d }
\def \fp{\int_M |u|^p\, \nu_g }
\def \mf{|u|^{p-2}}
\def \R{I\!\!R}
\def \dfp{\int_M |du|^p\, \d }

\def \fp{\int_M |u|^p\, \nu_g }

\def \sf{(S^n,\, can)}
\def \gf{\gamma_t^a \circ \phi }
 
\def \cun{(n+1)^{{p \over 2}-1}}
\def \cde{(n+1)^{1-{p \over 2}}}
\def \dis{\displaystyle}

\def \ba{\gamma_t^a}

\def \cqu{p \over 2}
\def \cci{p \over m}
\def \mo{|{p \over 2}-1|}
\def \ad{\sum_{i=1}^{n+1}}

\def \ep{\varepsilon}

\def \N{I\!\!N}

\begin{document}

\title{Conformal bounds for the first eigenvalue of the $p$-Laplacian}

\author{Ana-Maria Matei}
\address{Department of Mathematical Sciences, Loyola University New Orleans, USA; e-mail: amatei@loyno.edu}
\maketitle

\begin{abstract}: Let $M$ be a compact, connected, m-dimensional manifold without boundary and $p>1$. For $1<p\leq m$, we prove that the first eigenvalue $\la$ of the $p$-Laplacian is bounded on each conformal class of Riemannian metrics of volume one on $M$. For $p>m$, we show that any conformal class of Riemannian metrics on $M$ contains metrics of volume one with $\la$ arbitrarily large. As a consequence, we obtain that in two dimensions $\la$ is uniformly bounded on the space of Riemannian metrics of volume one if $1<p\leq 2$, respectively unbounded if $p>2$.
\end{abstract}

\vspace{0.4cm}

\vspace{0.3cm}

\bigskip

\section{Introduction}
\medskip

\noindent
Let $M$ be a compact m-dimensional manifold. All through this paper we will assume that $M$ is connected and without boundary. The $p$-Laplacian ($p>1$) associated to a Riemannian metric $g$ on $M$ is given by
$$
\Delta_{p} u=\delta(|du|^{p-2}du)\, ,
$$
\noindent
where $\delta = -div_g$ is the adjoint 
of $d$ for the $L^2$-norm induced by $g$ on the space of differential 
forms. This operator can be viewed as an extension of the Laplace-Beltrami operator which corresponds to $p=2$. The real numbers $\lambda$ for which the nonlinear partial differential equation  
$$
\Delta_{p} u =\lambda |u|^{p-2} u
$$
\noindent
 has nontrivial 
solutions are the {\em eigenvalues} of $\Delta_p$, and the associated solutions are the {\em 
eigenfunctions} of $\Delta_p$. Zero is an eigenvalue of $\Delta_p $, the associated eigenfunctions 
being the constant functions. The set of the nonzero eigenvalues is a nonempty, unbounded subset of 
$(0, \infty )$ \cite{GP}. The infimum $\la $ of this set is itself a positive eigenvalue, the {\em first eigenvalue} of $\Delta_p$, and has a Rayleigh type variational characterization \cite{Ve}:
$$
\la(M,g) = \inf \left \{ { {\dfp } \over  {\fp }} \; \mid \; u\in 
W^{1,p}(M)\setminus \{0\} \, , \; \s \mf u \d = 0\right \}\, ,
$$
\noindent
where $\nu_g$ denotes the Riemannian volume element associated to $g$.

The first eigenvalue of $\Delta_p$ can be viewed as a functional on the space of 
Riemannian metrics on M:
$$
g \mapsto \la(M,g)\, .
$$
\noindent
Since $\la$ is not invariant under dilatations ($\la(M,cg)=c^{-{p \over 2}}\la(M,g)$), a normalization is needed when studying the uniform boundedness of this functional. It is common to restrict $\la$ to  the set $\mathcal{M}(M)$ of Riemannian metrics of volume one on $M$. In the linear case $p=2$ this problem has been extensively studied in various degrees of generality. The functional $\lambda_{1,2}$ was shown to be uniformly bounded on $\mathcal{M}(M)$ in two dimensions \cite{her}, \cite{yaya}, \cite{LY}, and unbounded in three or more dimensions \cite{tan}, \cite{ura}, \cite{muto}, \cite{berber}, \cite{col}, \cite{dodziuk}. However, $\lambda_{1,2}$ becomes uniformly bounded when restricted to any conformal class of  Riemannian metrics in $\mathcal{M}(M)$ \cite{E-I1}. 

In the general case $p>1$, the functional $\la$ is unbounded on $\mathcal{M}(M)$ in three or more dimensions \cite{matei8}. In this paper we study the existence of uniform upper bounds for the restriction of $\la$ to conformal classes of Riemannian metrics in $\mathcal{M}(M)$:

$\bullet$ for $1<p\leq m$ we extend the results from the linear case and obtain an explicit upper bound for $\la$ in terms of $p$, the dimension $m$ and  the Li-Yau $n$-conformal volume. 

$\bullet$ for $p>m$, we consider first the case of the unit sphere $S^m$ and we construct Riemannian metrics in $\mathcal{M}(S^m)$, conformal to the standard metric $can$ and with $\la$ arbitrarily large. We use then the result on spheres to show that any conformal class of Riemannian metrics on $M$ contains metrics of volume one with $\la$ arbitrarily large. 

As a consequence, we obtain that in two dimensions, $\la$ is uniformly bounded on $\mathcal{M}(M)$ when $1<p\leq 2$, and unbounded when $p>2$.

\vspace{0.2cm} 

\section{The case $1<p\leq m\, :\,$ Li-Yau type upper bounds}
\medskip

\noindent
Let $g$ be a Riemannian metric on $M$ and denote by $[g]=\{fg\, | \, f\in C^{\infty}(M), \, f>0 \, \}$ 
the conformal class of $g$. Let $G(n)=\{\gamma \in \mbox{Diff}(S^n)\mid \gamma^*can \in [can]\}$  denote the group of conformal diffeomorphisms of $\sf $. 

For $n$ big enough, the Nash-Moser Theorem ensures (via the
stereographic projection) that  the set $I_n(M,[g])=\left\{\phi : M \to S^n\, \mid \, \phi^*can \in [g]  \right\}$ of conformal immersions from $ (M,g) $ to $\sf $ is nonempty.
The {\em n-conformal volume} of $[g]$ is defined by \cite{LY}: 
$$V_n^c(M,[g])= \, \inf_{\phi \in I_n(M,[g])} \sup_{\gamma \in G(n)}
{ Vol} \left (M, (\gamma \circ \phi )^{\ast}    can \right )\, ,$$
\noindent
where ${ Vol} \left (M, (\gamma \circ \phi )^{\ast}    can \right )$ denotes the volume of $M$ with respect to the induced metric $(\gamma \circ \phi )^{\ast}    can$. By convention, $V_n^c(M,[g])=\infty$ if $I_n(M,[g])=\emptyset$.
\noindent
\begin{theorem}\label{volconf}
Let $M$ be an $m$-dimensional compact manifold and $1<p\leq m$. For any metric $g \in \mathcal{M}(M)$ and any $n\in \N$ we have
 $$\la (M,g) \, \leq \, m^{\cqu }  (n+1)^{|{p \over 2}-1|}    V_n^c(M,[g])^{p \over 
 m}\,.
$$
\end{theorem}

\noindent
\begin{remark}
In the linear case $p=2$, this result was proved by Li and Yau 
\cite{LY} for surfaces and by El Soufi and Ilias \cite{E-I1} for 
higher dimensional manifolds.
\end{remark}
 
\begin{remark}
Theorem \ref{volconf} gives an explicit upper bound for $\la $, $1<p\leq m$, in the case of some particular   manifolds: the sphere $
{ S}^m$, the real projective space $ {\R \P}^m$, the complex projective 
space ${C \P}^d$, the equilateral torus
${{\T}_{eq}^2}$, the generalized Clifford torus  
$ {{ S}^r\left ( {\sqrt{r/{r+q}}}\right )\times { S}^q \left 
({\sqrt{q / 
{r+q}}}\right )}$, endowed with their canonical metrics. For these manifolds we 
have \cite{E-I1}: $V_n^c (M,[can])={Vol}(M, {{\lambda_{1,2}}\over m}can)$      
 for $n+1$ 
 greater or equal to the multiplicity of $\lambda_{1,2}$.
\end{remark}

\noindent
Using the relationships between the conformal volume and the genus of a 
compact surface \cite{E-I2} we obtain:

\noindent
\begin{corollary}\label{riemsurf}
Suppose $m=2$ and $1<p \leq 2$. Then for any metric ${g}\in \mathcal{M}(M)$
$$
\la (M,g)\leq k_p \left[{{genus(M)+3}\over {2}}\right]^{p\over 2}
\, ,$$
\noindent
where $[\,]$ denotes the integer part, $k_p=3^{|{p 
\over 2}-1|}(8\pi )^{p \over 2}$ if $M$ is orientable and $k_p=5^{|{p 
\over 2}-1|}(24\pi )^{p \over 2}$ if not. 
\end{corollary}

\begin{remark}
In the case $p=2$ and $M=S^2$, this result is the well known 
Hersch inequality \cite{her}. For higher genus surfaces, the upper bound 
of $\lambda_{1,2}$ in terms of the genus was obtained by El-Soufi and Ilias \cite{E-I2} by improving  a previous result of Yang and Yau \cite{yaya}.

\end{remark}

\noindent
In order to prove Theorem \ref{volconf} we need two Lemmas:
\noindent
\begin{lemma}\label{Degre}
Let $\phi     :(M,g) \to \sf$ be a smooth map whose level sets are of measure zero in $(M,g)$. 
Then for any $p>1$ there exists $\gamma \in G(n)$ such that     
$$\s |(\gamma \circ \phi)_i|^{p-2}    (\gamma \circ \phi)_i \d =0\, , \quad 1\leq i\leq n+1.
$$ 
\end{lemma}

\noindent
{\em Proof of Lemma \ref{Degre}}. Let $a \in S^n$ and denote by $ \pi_a$ the stereographic projection of pole $a$. Let $t \in (0,1] $  and $H_{{1-t} \over t}=e^{{1-t} \over t}\cdot Id_{\R^n}$ ( i.e. $H_{{1-t} \over t}$ is the
linear dilatation of $\R^n$ of factor $e^{{1-t} \over t}$). Let $\ba \in G(n)$, 
$\ba (x)=
\left\{
\begin{array}{lll}
 \pi_a^{-1} \circ H_{{1-t} \over t} \circ 
\pi_a (x) & \mbox{if} & x\in S^n \setminus \{ a \}\\
 a & \mbox{if} & x=a\, 
\end{array}
\right.
$ and consider the continuous map 
$$F   : (0,1] \times S^n \to \R^{n+1}$$ 
{\small $$
F(t,a)=  {1
\over {Vol(M,g)}} 
\left ( \s |(\gf )_1|^{p-2}(\gf )_1 \d ,\ldots ,\s |(\gf 
)_{n+1}|^{p-2}(\gf )_{n+1} \d \right ). 
$$}
\newline
For any  $x 
\in M \backslash \{ \phi^{-1}(-a)\} $ we have $\lim_{t \to 0^+} \gf(x)=a$. Since $\phi^{-1}(-a)$ is of measure 
zero in $M$, we can extend $F$ into a continuous function on $[0,1]\times S^n$ by setting
$${F(0,a)=\left (|a_1|^{p-2}a_1 ,\ldots 
,|a_{n+1}|^{p-2}a_{n+1}\right )}\,.
$$
\noindent
 The map $a \to F(0,a) $ is odd  on $S^n$, and since $\gamma_1^a =Id_{S^n}$, the map $a \to F(1,a)$ 
is constant. 
Assume $|| F(t,a) ||  \neq 0$ for any $(t,a)\in [0,1]\times S^n$. Then the map 
$$G:[ 0,1 ] \times S^n \to S^n$$ 
$$
G(t,a)={{F(t,a)}\over || F(t,a) || }
$$
\noindent
gives a homotopy between the odd map $a \to G(0,a)$
and the constant map $a \to G(1,a) $, and this is impossible.
Hence there exists $(t,a)\in [0,1] \times S^n $ such that 
$||F(t,a)||=0$, i.e.    
$
\s | \left ( \gf \right )_i|^{p-2}\left ( \gf\right )_i \d =0,\; \; \; 
1\leq i\leq n+1$.
\hbox{}\hfill$\Box$
\vspace{0.2cm}

\noindent
\begin{lemma}\label{p-energy}
Suppose $g \in \mathcal{M}(M)$ and let $\phi     :(M,g) \to \sf$ be a smooth map whose level sets are of measure zero in $(M,g)$. Then there exists $\gamma \in G(n)$ such that  
$$ 
\la (M,g) \, \leq \,  {(n+1)^{\mo }} \int_M |d(\gamma \circ \phi)|^p\nu_g  
\, ,$$
where $|d(\gamma \circ \phi)|$ denotes the Hilbert-Schmidt norm of  $d(\gamma \circ \phi)$.
\end{lemma}
\vspace{0.2cm}

\noindent 
{\em Proof of Lemma \ref{p-energy}}. Lemma \ref{Degre} implies there exists $\gamma \in G(n)$ such that $\psi    =\gamma \circ 
\phi : M\to S^n$ verifies  $ {\s |\psi_i|^{p-2}   \psi    _i\d=0}$, 
$1\leq i\leq n+1$. The variational characterization for $\la(M,g)$ implies that 
$\la(M,g)\leq \dis {{\int_M |d\psi_i|^p\nu_g}\over {\int_M|\psi_i|^p\nu_g}}$, $1\leq i\leq n+1$. Then
\begin{equation}\label{var-psi}
\la(M,g)\leq \dis {\int_M {\sum_{i=1}^{n+1} |d\psi_i|^p\,\nu_g}\over {\int_M{\sum_{i=1}^{n+1}|\psi_i|^p\,\nu_g}}}\, .
\end{equation}
\newline
$\bullet$ {\em Case 1:} ${p \geq 2}.$ It is straightforward that
\begin{equation}\label{num1}
\ad | d\psi_i |^p = \ad (|d\psi_i |^2)^{\cqu }\leq \left (\ad |d\psi_i|^2  
\right )^{\cqu }=|d\psi |^p \, . 
\end{equation} 
\noindent
On the other hand
\begin{equation}\label{den1}
\sum_{i=1}^{n+1}|\psi_i|^p\geq (n+1)^{1-{p\over 2}}\left(\sum_{i=1}^{n+1}|\psi_i|^2\right)^{p\over 2}=
 (n+1)^{1-{p\over 2}}\, ,
\end{equation}
\noindent
where we have used the fact that $x\to x^{p\over 2}$ is convex and that $\sum_{i=1}^{n+1}|\psi_i|^2=1$. Replacing (\ref{num1}) and (\ref{den1}) in (\ref{var-psi}) we obtain
$$
\la (M,g)  \leq  \cun \s |d\psi |^p \d  \, .
$$ 
\noindent
$\bullet$ {\em Case 2:} $1< p<2$. Since $|\psi_i|\leq 1$ we have $|\psi_i|^2 \leq |\psi_i |^p$ and
\begin{equation}\label{den2}
1=Vol(M,g) =\s \ad |\psi_i|^2 \d  \leq \s \ad |\psi_i|^p  \d 
\end{equation}
\noindent
On the other hand 
\begin{equation}\label{num2}
\ad |d\psi_i |^p = \ad (| d\psi_i |^2)^{\cqu }
\leq \cde \left (\ad | d\psi_i |^2 \right )^{\cqu }=\cde | d\psi |^p,
\end{equation}
\noindent
where the inequality follows from the concavity of $x\to x^{p\over 2}$. Replacing (\ref{den2}) and (\ref{num2}) in (\ref{var-psi}) we obtain
$$
\la (M,g) \leq \cde \int_{M}|d\psi|^p\, \nu_{g}\, .
$$
\hbox{}\hfill$\Box$
\vspace{0.2cm}

\noindent
{\em Proof of Theorem \ref{volconf}}. Let  $\phi :(M,g) \to (S^{n},can )$ be a conformal immersion. From Lemma \ref{p-energy} we have that there exists $\gamma \in G(n)$ such that
$$
\la (M,g)  \leq (n+1)^{|{p \over 2}-1|} \s |d (\gamma \circ \phi)|^p \nu_{ g}\, .
$$
\newline
Since $g\in \mathcal{M}(M)$, H\"older's inequality implies    
$$\s |d (\gamma \circ \phi) |^p \nu_{g} \, \leq \, \left ( \s | d (\gamma \circ \phi)|^m \nu_{ g}\right )^{\cci }.$$
\noindent
On the other hand since $ \gamma \circ \phi :(M, g)\to (S^n,can)$ is a conformal immersion, $ (\gamma \circ \phi)^*can= 
{{|d (\gamma \circ \phi)|^2}\over m} g$ and we have
 $$ \s |d  (\gamma \circ \phi)|^m \nu_{ g} =m^{m \over 2}{Vol}(M,  (\gamma \circ \phi)^\ast can)\, \leq \,
m^{m \over 2}\sup_{\gamma \in G(n)}{Vol}(M,(\gamma \circ \phi)^{\ast }can ).$$
\noindent
Combining the inequalities above we obtain:
$$
\la(M, g)\leq \,m^{p \over 2}(n+1)^{|{p \over 2}-1|}
\left(\sup_{\gamma \in G(n)}{Vol}(M,(\gamma \circ \phi)^{\ast }can )\right)^{p\over m}
$$
Taking the infimum over all $\phi \in I_n(M,[g])$ we obtain the desired inequality.
\hbox{}\hfill$\Box$
\vspace{0.2cm}

\noindent
{\em Proof of Corollary \ref{riemsurf}}. In the case of surfaces, the $n$-conformal volume  is bounded above by a constant depending only on the genus of the surface \cite{E-I2}. If $M$ is orientable we have 
$$
V_n^c(M,[g]) \leq 4\pi   \left[{{genus(M)+3}\over {2}}\right]   
\quad \mbox{for} \quad n \geq 2\, .
$$
\noindent
If $M$ is non orientable,
$$
V_n^c(M,[g]) \leq 12\pi \left[{{genus (M)+3}\over{2}}\right] \quad \mbox{for}\quad n \geq 4 \, .
$$
\noindent
Theorem \ref{volconf} implies now the desired result with $k_p=3^{|{p 
\over 2}-1|}(8\pi )^{p \over 2}$ when $M$ is orientable and $k_p=5^{|{p 
\over 2}-1|}(24\pi )^{p \over 2}$ when $M$ is non orientable.
\hbox{}\hfill$\Box$

\section{The case $p> m$}
\medskip

\noindent
For the sake of self-containedness we include here the variational characterizations for the first eigenvalues for the Dirichlet and the Neumann problems for $\Delta_p$. Let $\Omega$ be a domain in $M$ and consider the Dirichlet problem:
$$
\left \{ \begin{array}{ll}
{\Delta_pu} &  {= \lambda \, |u|^{p-2} u \quad {\rm in} \; \Omega} \\
{u} &  {=0\quad \quad {\mbox {on}}\; { \partial \Omega }} \, .
\end{array}
\right.
$$
The infimum $\lambda_{1,p}^D(\Omega,g)$ of the set of eigenvalues for this problem is itself a positive eigenvalue with the variational characterization
$$
\la^D(\Omega,g)=\inf \left \{ { {\int_{\Omega} |du|^p \nu_g}
\over  {\int_{\Omega}|u|^p \nu_g }} \; | \; u \in 
W_0^{1,p}(\Omega)\setminus\{0\} \,\right \}.
$$
Consider now the Neumann problem on $\Omega$:
$$
\left \{ \begin{array}{ll}
\dis {\Delta_pf} &  {= |f|^{p-2} f \quad {\rm in} \; \Omega} \\
df(\eta) & {=0\quad \quad {\mbox {on}}\; { \partial \Omega }} \, ,
\end{array}
\right.
$$
\noindent
where $\eta$ denotes the exterior unit normal vector field to $\partial \Omega$. Here too, the infimum $\la^N(\Omega,g)$ of the set of nonzero eigenvalues is a positive eigenvalue with the variational characterization
$$
\la^N(\Omega,g):=\inf \left \{ {{\int_{\Omega}|df|^p \, \nu_g } \over 
 {\int_{\Omega}|f|^p \, \nu_g }} \mid f \in 
W^{1,p}({\Omega},g)\setminus\{0\} \;, \; \int_{\Omega}|f|^{p-2}f \, \nu_g = 0\right \}.
$$
\vspace{0.2cm}

\noindent
We consider first the case of $(S^m, [can])$:

\begin{theorem}\label{boundedness-sphere}
For any $p>m$, $S^m$ carries Riemannian metrics of
volume one, conformal to the standard metric $can$, with $\la$ arbitrarily large.
\end{theorem}

\noindent
{\em Proof of Theorem \ref{boundedness-sphere}}. Let $r \in
[0,{\pi}]$, denote the geodesic distance on $(S^m, can)$ w.r.t. a point $x_0\in S^m$. Let $\ep >0$ and define a radial function $f_{\ep}:S^m \to \R$ by
\begin{equation}\label{singular}
f_{\ep}(r)={{\ep}^{{4p}\over{m(p-m)}}} \cdot \chi_{[0, {{\pi }\over 2}-\ep]\cup [ {{\pi }\over 2}+\ep, \pi]}(r) +  \chi_{({{\pi}\over 2}-\ep, {{\pi}\over 2}+\ep)}(r)\,.
\end{equation}
\noindent
Let 
\begin{multline*}
\la(\ep)=\inf \left \{ R_{\ep}(u):= {
{\int_{S^m}|du|^p\,f_{\ep}^{{m-p}\over 2}  \nu_{can} }
\over {\int_{S^m}|u|^p \, f_{\ep }^{m\over 2}\nu_{can}}} \mid 
u\in W^{1,p}(S^m)\setminus\{0\} \, , \right.
 \\
\left. \int_{S^m}|u|^{p-2}u \, |f_{\ep }|^{m\over 2}\nu_{can}=0
 \right\} .
\end{multline*}
\noindent
We will show first that
\begin{equation}\label{lim}
\limsup_{\ep \to 0}\la(\ep)\cdot \ep^{p\over m}=\infty\, .
\end{equation}
\noindent
Classical density arguments imply that there exists 
$u_{\ep}\in W^{1,p}(S^m)\setminus\{0\}$ with $\int_{S^m}|u_{\ep}|^{p-2}u_{\ep}f_{\ep }^{m\over 2}\nu_{can}=0$ such that $\la(\ep)=R_{\ep}(u_{\ep})$.  Let $\bar{u}_{\ep}:S^m\to \R$ be a radial function defined by
\begin{equation}\label{denom}
\bar{u}_{\ep}^p(r)={1\over {V}}\int_{S^{m-1}}|u_{\ep}(r,\cdot)|^p \nu_{can}
\end{equation}
\noindent
where $V=Vol(S^{m-1}, can)$. Differentiating w.r.t. $r$ we obtain
$$
p\bar{u}_{\ep}^{p-1}\bar{u}_{\ep}'={p\over {V}} \int_{S^{m-1}}|u_{\ep}|^{p-2}u_{\ep}
{{\partial u_{\ep}}\over {\partial r}} \nu_{can}\, .
$$
\noindent
By H\"older's inequality we obtain
$$
\bar{u}_{\ep}^{p-1}|\bar{u}_{\ep}'|\leq {1\over V} \int_{S^{m-1}}|u_{\ep}|^{p-1}\left |{{\partial u_{\ep}}\over {\partial r}} \right | \nu_{can}
\leq {1\over V}
\left (
 \int_{S^{m-1}}|u_{\ep}|^{p}\nu_{can}
\right )^{{p-1}\over p}
\cdot
\left (
 \int_{S^{m-1}}\left | {{\partial u_{\ep}}\over {\partial r}}  \right |^{p}\nu_{can}
\right )^{1\over p}.
$$
\noindent
It follows that
\begin{equation}\label{num}
|\bar{u}_{\ep}'|^p \leq {1\over V}\int_{S^{m-1}}\left |{{\partial u_{\ep}}\over {\partial r}} \right |^{p}\nu_{can}\leq {1\over V}\int_{S^{m-1}} |du_{\ep} |^{p}\nu_{can}
\end{equation}
\noindent
On the other hand
$$
\begin{array}{ll}
\displaystyle{\int_{S^m}|\bar{u}_{\ep}|^p \, f_{\ep }^{m\over 2}\nu_{can}}& =
\displaystyle{V\cdot \int_0^{\pi}|\bar{u}_{\ep}|^p \, f_{\ep }^{{m}\over 2}\sin r^{m-1}dr}\\ 
& \\
& =
\displaystyle{\int_0^{\pi}\left[\int_{S^{m-1}}|u_{\ep}|^p \, \nu_{can}\right]f_{\ep }^{m\over 2}\sin r ^{m-1}dr}\\
 & \\
 & =
\displaystyle{\int_{S^m}|u_{\ep}|^p \, f_{\ep }^{m\over 2}\nu_{can}},
\end{array}
$$
\noindent
where the second equality follows from (\ref{denom}). Similarly (\ref{num}) implies
$$
\int_{S^m}|\bar{u}_{\ep}'|^p \, f_{\ep }^{{m-p}\over 2}\nu_{can}
\leq \int_{S^m}|du_{\ep}|^p \, f_{\ep }^{{m-p}\over 2}\nu_{can}\, .
$$
\noindent
In particular, we obtain that 
$\bar{u}_{\ep}\in W^{1,p}(S^m)$ and 
$$
\begin{array}{ll}
\la(\ep)=R_{\ep}(u_{\ep})& \geq \displaystyle{{\int_{S^m}|\bar{u}_{\ep}'|^p \, f_{\ep }^{{m-p}\over 2}\nu_{can}}\over 
{\int_{S^m}|\bar{u}_{\ep}|^p \, f_{\ep }^{m\over 2}\nu_{can}}}\\ 
& \\
& \geq \min \left \{ 
\displaystyle{{\int_{S_+^m}|\bar{u}_{\ep}'|^p \, f_{\ep }^{{m-p}\over 2}\nu_{can}}\over 
{\int_{S_+^m}|\bar{u}_{\ep}|^p \, f_{\ep }^{m\over 2}\nu_{can}}}\, , \, 
\displaystyle{{\int_{S_-^m}|\bar{u}_{\ep}'|^p \, f_{\ep }^{{m-p}\over 2}\nu_{can}}\over 
{\int_{S_-^m}|\bar{u}_{\ep}|^p \, f_{\ep }^{m\over 2}\nu_{can}}} 
\right \}\, ,
\end{array}
$$
\noindent
where $S_+^m, S_{-}^m$ denote the hemispheres centered at $x_0$, respectively $-x_0$. Without loss of generality we may assume that
\begin{equation}\label{baru}
\la(\ep)\geq {{\int_{S_+^m}|\bar{u}_{\ep}'|^p \, f_{\ep }^{{m-p}\over 2}\nu_{can}}\over 
{\int_{S_+^m}|\bar{u}_{\ep}|^p \, f_{\ep }^{m\over 2}\nu_{can}}}\, .
\end{equation}
\noindent
Let $w_{\ep}\in W^{1,p}(S^m_+)$, 
$w_{\ep}=
\left\{
\begin{array}{lll}
\bar{u}_{\ep} & \mbox{on}  & [0, {{\pi}\over 2}-\ep]\\ 
\bar{u}_{\ep}({{\pi}\over 2}-\ep) & \mbox{on} & ({{\pi}\over 2}-\ep, {{\pi}\over 2}]
\end{array}
\right.$ and $v_{\ep}=\bar{u}_{\ep}-w_{\ep}$. Then $v_{\ep}=0$ on $[0,{{\pi}\over 2}-\ep]$ and $w_{\ep}'=0$ on $({{\pi}\over 2}-\ep, {{\pi}\over 2})$. Since $v_{\ep}'$ and $w_{\ep}'$ have disjoint supports, we have $|\bar{u}_{\ep}'|^p=  |v_{\ep}'|^p+|w_{\ep}'|^p$. On the other hand
$|\bar{u}_{\ep}|^p=|v_{\ep}+w_{\ep}|^p\leq 2^{{p}-1}( |v_{\ep}|^p+|w_{\ep}|^p)$. Then (\ref{baru}) and (\ref{singular}) imply
$$
\begin{array}{ll}
\la(\ep) & \geq  2^{1-{p}}\displaystyle{{\int_{S_+^m} (|v_{\ep}'|^p+|w_{\ep}'|^p ) f_{\ep}^{{m-p}\over 2}\nu_{can}}\over {{\int_{S_+^m} (|v_{\ep}|^p
+|w_{\ep}|^p) f_{\ep}^{{m}\over 2}\nu_{can}}}}\\
 & \\
  & =2^{1-{p}}\displaystyle{{\int_{S^m_+} |v_{\ep}'|^p \nu_{can}+ 
  {\ep}^{-{{2p}\over m}}\int_{S^m_+}|w_{\ep}'|^p \nu_{can}}\over {{\int_{S^m_+} |v_{\ep}|^p \nu_{can}+ \int_{S^m_+}|w_{\ep}|^p f_{\ep}^{{m}\over 2}\nu_{can}}}}
\end{array}
$$
Quite to multiply $\bar{u}_{\ep}$ by a constant we may assume ${\int_{S^m_+} |v_{\ep}|^p \nu_{can}+ \int_{S^m_+}|w_{\ep}|^p f_{\ep}^{{m}\over 2}\nu_{can}}=1$ and the inequality above becomes
\begin{equation}\label{norm}
\la(\ep)\geq 2^{1-{p}}{\int_{S^m_+} |v_{\ep}'|^p \nu_{can}+ 
  {\ep}^{-{{2p}\over m}}\int_{S^m_+}|w_{\ep}'|^p \nu_{can}}
\end{equation}
\noindent
$\bullet$ {\em Case 1}:  $\limsup_{\ep \to 0}\int_{S^m_+}|w_{\ep}'|^p \nu_{can}>0$. 
\newline
Inequality (\ref{norm}) implies that $\la(\ep)\ep^{p\over m}\geq 2^{1-{p}}
{\ep}^{-{{p}\over m}}\int_{S^m_+}|w_{\ep}'|^p \nu_{can}$, and therefore (\ref{lim}) is verified.
\newline
$\bullet$ {\em Case 2}: $\lim_{\ep \to 0}\int_{S^m_+}|w_{\ep}'|^p \nu_{can}=0$. 
\newline
Then we may find a sequence $\ep_n\to 0$ such that $w_{\ep_n}\to c$ strongly in 
$L^{p}(M)$, where $c$ is a constant. In particular since $p>m$, $\{f_{\ep_n}\}$ is uniformly bounded and we have
$\lim_{n \to \infty}\int_{S^m_+} f_{\ep_n}^{{m}\over 2}\nu_{can}=0$. It follows that 
$\lim_{n \to \infty}\int_{S^m_+}|w_{\ep_n}|^p f_{\ep_n}^{{m}\over 2}\nu_{can}=\lim_{n\to \infty}\int_{S^m_+}(|w_{\ep_n}|^p-|c|^p) f_{\ep_n}^{{m}\over 2}\nu_{can}+
|c|^p\lim_{n\to \infty}\int_{S^m_+} f_{\ep_n}^{{m}\over 2}\nu_{can}=0$. Hence for $\ep_n $ small enough, ${\int_{S^m_+} |v_{\ep_n}|^p \nu_{can}}=1-\int_{S^m_+}|w_{\ep_n}|^p f_{\ep}^{{m}\over 2}\nu_{can}\geq {1\over 2}
$ and (\ref{norm}) implies 
\begin{equation}\label{sin}
\begin{array}{ll}
\la(\ep_n) &\geq 2^{1-{p\over 2}}\int_{S^m_+} |v_{\ep_n}'|^p \nu_{can}\geq 2^{-{p\over 2}}\displaystyle{{\int_{S^m_+} |v_{\ep_n}'|^p \nu_{can}}\over {{\int_{S^m_+} |v_{\ep_n}|^p \nu_{can}}}}
=2^{-{p\over 2}}
\displaystyle{{\int_{{{\pi}\over 2}-\ep_n}^{{\pi}\over 2} |v_{\ep_n}'|^p \sin r ^{m-1} dr}\over {\int_{{{\pi}\over 2}-\ep_n}^{{\pi}\over 2} |v_{\ep_n}|^p \sin r ^{m-1}dr}}\\
& \\ 
& \geq 2^{-{p\over 2}}[\sin({{ {\pi}\over 2}-\ep_n})]^{m-1}\displaystyle{{\int_{{{\pi}\over 2}-\ep_n}^{{\pi}\over 2} |v_{\ep_n}'|^p  dr}\over {\int_{{{\pi}\over 2}-\ep_n}^{{\pi}\over 2} |v_{\ep_n}|^p dr}}.
\end{array}
\end{equation}
\noindent
Let $\bar{v}_{\ep_n} \in W_0^{1,p}(-\ep_n, \ep_n)$ be an even function such that $\bar{v}_{\ep_n}(s)=v_{\ep_n}(s+{{\pi}\over 2}-\ep_n)$ for $0\leq s\leq {\ep_n}$. We have then
\begin{equation}\label{dir}
\displaystyle{{\int_{{{\pi}\over 2}-\ep_n}^{{\pi}\over 2} |v_{\ep_n}'|^p  dr}\over {\int_{{{\pi}\over 2}-\ep_n}^{{\pi}\over 2} |v_{\ep_n}|^p dr}}= \displaystyle{{\int_{0}^{\ep_n} |\bar{v}_{\ep_n}'|^p  dr}\over {\int_{0}^{\ep_n} |\bar{v}_{\ep_n}|^p dr}}=\displaystyle{{\int_{-\ep_n}^{\ep_n} |\bar{v}_{\ep_n}'|^p  dr}\over {\int_{-\ep_n}^{\ep_n} |\bar{v}_{\ep_n}|^p dr}}\geq \la^D(-\ep_n, \ep_n)=\ep_n^{-p}\la^D(-1,1)\, .
\end{equation}
\noindent
Inequalities (\ref{sin}), (\ref{dir}) imply $\la(\ep_n)\geq \ep_n^{-p}\la^D(-1,1)$ and (\ref{lim}) is verified again.
\vspace{0.1cm}

\noindent
Fix now $\ep>0$ and let $\tilde{f}_{\ep}\in C^{\infty}(S^m)$, radial with respect to $x_0$ and such that $\tilde{f}_{\ep}\leq {f}_{\ep}$, $\tilde{f}_{\ep}(r)={f}_{\ep}(r)=1$ on $[{{\pi}\over 2}-{{\ep}\over 2}, {{\pi}\over2}+{{\ep}\over 2}]$ and $\tilde{f}_{\ep}{(\pi-r)}=\tilde{f}(r)$. Then
\begin{equation}\label{vol}
\begin{array}{ll}
Vol (S^m, \tilde{f}_{\ep} can)& ={\int_{S^m} \tilde{f}_{\ep}^{m\over 2}\nu_{can}}=
{\int_{S^{m-1}}
\int_{-{{\pi}\over 2}}^{{\pi}\over 2}\tilde{f}_{\ep}^{m\over 2}\sin r^{m-1}\, dr\, \nu_{can}}\\
& \\
 & > V{\int_{{{\pi}\over 2}-{{\ep}\over 2}}^{{{\pi}\over 2}+{{\ep}\over 2}}\sin r^{m-1}dr}\\
 & \\
& >\ep V[\sin({{\pi}\over 2}-\ep)]^{m-1},\quad \mbox{where}\;V=Vol(S^{m-1},can).
\end{array}
\end{equation}
We will compare now $\la(S^m, \tilde{f}_{\ep}can)$ and $\la(\ep)$. Let $\tilde{u}_{\ep}$ be an eigenfunction for $\la(S^m, \tilde{f}_{\ep}can)$ and denote by $\tilde{u}_{\ep}^+, \tilde{u}_{\ep}^-$ the positive, respectively, the negative part of $\tilde{u}_{\ep}$.  Then \cite{matei1}
 $$
\la(S^m, \tilde{f}_{\ep}can)={
{\int_{S^m}|d\tilde{u}_{\ep}^+|^p\,\tilde{f}_{\ep}^{{m-p}\over 2}  \nu_{can} }
\over {{\int_{S^m}}|\tilde{u}_{\ep}^+|^p \, \tilde{f}_{\ep }^{m\over 2}\nu_{can}}}
=
{
{\int_{S^m}|d\tilde{u}_{\ep}^-|^p\,\tilde{f}_{\ep}^{{m-p}\over 2}  \nu_{can} }
\over {{\int_{S^m}}|\tilde{u}_{\ep}^-|^p \, \tilde{f}_{\ep }^{m\over 2}\nu_{can}}}
$$
\noindent
Let $t \in \R$ and $\tilde{u}_{\ep,t}=t
\tilde{u}_{\ep}^++\tilde{u}_{\ep}^-$. Then there is $t_0$ such that $\int_{S^m}|\tilde{u}_{\ep,t_0}|^{p-2}\tilde{u}_{\ep,t_0}f_{\ep }^{m\over 2}\nu_{can}=0$ and the equation above implies
\begin{equation}\label{compare}
\la(S^m, \tilde{f}_{\ep}can)={
{\int_{S^m}|d\tilde{u}_{\ep,t_0}|^p\,\tilde{f}_{\ep}^{{m-p}\over 2}  \nu_{can} }
\over {{\int_{S^m}}|\tilde{u}_{\ep,t_0}|^p \, \tilde{f}_{\ep }^{m\over 2}\nu_{can}}}
\geq {
{\int_{S^m}|d\tilde{u}_{\ep,t_0}|^p\,{f}_{\ep}^{{m-p}\over 2}  \nu_{can} }
\over {{\int_{S^m}}|\tilde{u}_{\ep,t_0}|^p \, {f}_{\ep }^{m\over 2}\nu_{can}}}
\geq 
\la(\ep)\, ,
\end{equation}
where the first inequality follows from the fact that $\tilde{f}_{\ep}\leq f_{\ep}$ and the second from the variational characterization for $\la(\ep)$.
Inequalities (\ref{vol}),  (\ref{compare}) and (\ref{lim}) yield
$$
\limsup_{\ep \to 0}\la(S^m, \tilde{f}_{\ep} can)Vol(S^m,\tilde{f}_{\ep}can)^{p\over m} \\ \geq 
V^{p\over m}\cdot \limsup_{\ep \to 0}\la({\ep})\cdot \ep^{p\over m}=\infty .
$$
\noindent
Finally, let $h_{\ep}=Vol(S^m, \tilde{f}_{\ep}can)^{-{2\over m}}\tilde{f}_{\ep}$. We have then
\begin{equation*}
\hspace{1cm}Vol(S^m, h_{\ep}can)=1\quad \mbox{and}\quad 
\limsup_{\ep \to 0}\la(S^m, h_{\ep} can)=\infty\, .
\hspace{2cm}\Box
\end{equation*}

\noindent
We will extend the construction from $(S^m,[can])$ to $(M,[g])$ by means of the first eigenvalue for the Neumann problem for $\Delta_p$ on a domain $\Omega$ in $M$.

\begin{theorem}\label{boundedness-general}
Let $(M, g)$ be a compact Riemannian manifold of dimension $m$. Then for any $p>m$, $[g]$ contains Riemannian metrics of volume one with $\la$ arbitrarily large.
\end{theorem}

\noindent
{\em Proof of Theorem \ref{boundedness-general}}. Let $r$ denote the geodesic distance on $(S^m,can)$ w.r.t. a point $x_0$. Let $f\in C^{\infty}(S^m)$ be a function  radial w.r.t.  $x_0$, such that  $f(r)=f(\pi-r)$ and $Vol(S^m,fcan)=1$. As before, let $S_+^m$ denote the hemisphere centered at $x_0$. Let $v$ be an eigenfunction for 
$\la^N(S_+^m, {f}can)$ and let $w\in W^{1,p}(S^m)$, $
w(r)=
\left\{
\begin{array}{lll}
v(r)& \mbox{if}& 0\leq r\leq {{\pi}\over 2}\\
v(\pi -r) & \mbox{if} & {{\pi}\over 2}<r\leq \pi
\end{array}
\right.
$. Then $\int_{S^m}  |w|^{p-2}w {f}^{m\over 2}\nu_{can}=2 \int_{S_+^m}  |v|^{p-2}v {f}^{m\over 2}\nu_{can}=0$ and the variational characterization for $\la(S^m, {f}can)$ implies
\begin{equation}\label{new-sph}
\la(S^m, {f}can)\leq {{\int_{S^m}|dw|^p{f}^{{m-p}\over 2}\nu_{can}}\over
 {\int_{S^m}|w|^p{f}^{{m}\over 2}\nu_{can}}}={{\int_{S_+^m}|dv|^p{f}^{{m-p}\over 2}\nu_{can}}\over
 {\int_{S_+^m}|v|^p{f}^{{m}\over 2}\nu_{can}}}=\la^N(S_+^m, {f\,}can)
\end{equation}
\noindent
Let  $\Omega$ be a domain in $M$ such that there exists a diffeomorphism $\Phi: \Omega\to S_+^m$. We may assume $\Omega$ is included in the open region of a local chart of $M$. In this chart we have 
$\nu_g=\sqrt{\det(g_{ij})}dx^1\wedge dx^2\wedge\ldots \wedge dx^m$ and $
\nu_{\Phi^* can}=\sqrt{\det((\Phi^*can)_{ij})}dx^1\wedge dx^2\wedge\ldots \wedge dx^m$.
There exist positive constants $c_1,c_2$ such that 
\begin{equation}\label{const}
c_1\sqrt{\det(g_{ij})}\leq \sqrt{\det((\Phi^*can)_{ij})} \leq c_2 \sqrt{\det(g_{ij})}\quad \mbox {on}\; \Omega\,.
\end{equation}
\noindent
We will compare now $\la^N(S_+^m, {f\,}can)$ and $\la^N (\Omega,(f\circ \Phi)g)$. Note first that 
since $\Phi$ is an isometry between $(\Omega,(f\circ \Phi)\Phi^*can)$ and $(S^m_+, fcan)$ we have
\begin{equation}\label{new}
\la^N(S^m_+, fcan)=\la^N (\Omega,(f\circ \Phi)\Phi^*can)
\end{equation}
\noindent
Let $u$ be an eigenfunction for $\la^N(\Omega, (f\circ \Phi)g)$ and denote by $u^+$, $u^-$ the positive, respectively, the negative part of $u$. Then there is $s\in \R$ such that the function $u_s=su^++u^-$ verifies $\int_{\Omega}|u_s|^{p-2}u_s (f\circ \Phi)^{m\over 2} \nu_{\Phi^*can} =0$. Furthermore
\begin{equation}\label{compare-new}
\begin{array}{ll}
\la^N(\Omega, (f\circ \Phi)g) &
  =\displaystyle{{\int_{\Omega}|du_s|^p{(f\circ \Phi)}^{{m-p}\over 2}\nu_{g}}\over
 {\int_{\Omega}|u_s|^p{(f\circ \Phi)}^{{m}\over 2}\nu_{g}}}\geq  {{c_1}\over {c_2}}\displaystyle{{\int_{\Omega}|du_s|^p{(f\circ \Phi)}^{{m-p}\over 2}\nu_{\Phi^*can}}\over
 {\int_{\Omega}|u_s|^p{(f\circ \Phi)}^{{m}\over 2}\nu_{\Phi^*can}}}\\
 & \\
 & \geq \displaystyle{{c_1}\over {c_2}}\la^N(\Omega, (f\circ \Phi)\Phi^*can)\, ,
 \end{array}
\end{equation}
\noindent
where the first inequality follows from (\ref{const}) and the second from the variational characterization of $\la^N(\Omega, (f\circ \Phi)\Phi^*can)$.
From (\ref{new-sph}), (\ref{new}) and (\ref{compare-new}) we obtain
\begin{equation}\label{compare-mani}
\la^N(\Omega, (f\circ \Phi)g)\geq {{c_1}\over {c_2}}\la(S^m, fcan)\, .
\end{equation}
Let now $\delta >0$; there is an extension $\widetilde {f\circ \Phi}$ of ${f\circ \Phi}$ on the entire manifold $M$ such that the metric $\tilde g =\widetilde {f\circ \Phi} g$ verifies \cite{matei4}: $
{\la (M,\tilde g)>\la^N(\Omega,(f\circ \Phi)g)- \delta}$. Inequality (\ref{compare-mani}) implies 
\begin{equation}\label{msfd}
\la (M,\tilde g)>{{c_1}\over {c_2}}\la(S^m, fcan)-\delta
\end{equation}
\noindent
On the other hand
\begin{multline}\label{volume}
Vol(M,\tilde g)> Vol(\Omega,(f\circ \Phi)g) \geq {1\over{c_2}}Vol (\Omega, (f\circ \Phi)\Phi^*can)
\\ ={1\over{c_2}}Vol (S^m_+,fcan)={1\over{2c_2}}Vol (S^m,fcan)={1\over{2c_2}}.
\end{multline}
\noindent
Let $K>0$; from the proof of Theorem \ref{boundedness-sphere}  we may assume that  $f$ is chosen  such that  $\la(S^m, {f}can)>2^{{p\over m}+1}{{c_1^{-1}}{c_2^{{p\over m}+1}}}K$.
For $\delta$ small enough such that $(2c_2)^{-{p\over m}}\delta<K$, inequalities (\ref{msfd}) and (\ref{volume}) imply
$$
\begin{array}{ll}
\la(M,\tilde g)Vol(M, \tilde g)^{p\over m}&\geq
[{{c_1}\over {c_2}}\la(S^m, fcan)-\delta)](2c_2)^{-{p\over m}}>K\, .
\end{array}
$$
\noindent
Finally, let $h={Vol(M,\tilde{g})^{-{2\over m}}}\tilde{g}$. Then $h \in [g]$, $Vol(M, h)=1$ and $\la(M,h)>K$.
\hbox{}\hfill$\Box$

\end{document}